\documentclass{elsart}

\usepackage{amssymb,amsmath}
\theoremstyle{plain}

\newtheorem{exe}{Example}
\newtheorem{defi}{Definition}

\theoremstyle{definition}
\newcommand{\mem}[1]{\in\mathbb{#1}}

\begin{document}

\begin{frontmatter}



\title{A Geometric Proof to Cantor's Theorem and an Irrationality Measure for Some Cantor's Series}

\author{Diego Marques}
\address{Departamento de Matem\'{a}tica, Universidade Federal do Cear\'{a}, Fortaleza, Cear\'{a}, Brazil}
\ead{diego@mat.ufc.br}

\begin{abstract}
Generalizing a geometric idea due to J. Sondow, we give a geometric proof for the Cantor's Theorem. Moreover, it is given an irrationality measure for some Cantor series. 
\end{abstract}

\end{frontmatter}

\section{Introduction} 
In $2006$, Jonathan Sondow gave a nice geometric proof that $e$ is irrational. Moreover, he said that a generalization of his construction may be used to prove the Cantor's theorem. But, he did not do that in his paper, see \cite{Sondow}. So we give a geometric proof to Cantor's theorem using a generalization to Sondow's construction. After, it is given an irrationality measure for some Cantor series, for that we generalize the Smarandache function. Also we give an irrationality measure for $e$ that is a bit better than the given one in \cite{Sondow}.

\section{Cantor's Theorem} 

\begin{defi}
Let $a_0,a_1,...,b_1,b_2,...$ be sequences of integers that satisfy the inequalities $b_n \geq 2$, and $0 \leq a_n \leq b_n - 1$ if $n \geq 1$. Then the convergent series
\begin{equation}\label{cantor}
\theta := a_0 + \displaystyle\frac{a_1}{b_1} + \displaystyle\frac{a_2}{b_1b_2} + \displaystyle\frac{a_3}{b_1b_2b_3} + \ldots.
\end{equation}
is called $\textit{Cantor series}.$
\end{defi}
\begin{exe}
The number $e$ is a Cantor series. For see that, take $a_0 = 2, a_n = 1, b_n = n + 1$ for $n \geq 1.$
\end{exe}

We recall the following theorem due to Cantor \cite{Cantor}.

\textbf{Theorem (Cantor)} \textit{Let $\theta$ be a Cantor series. Suppose that each prime divides infinitely many of the $b_n$. Then $\theta$ is irrational if and only if both $a_n > 0$ and $a_n < b_n - 1$ hold infinitely often.}

\textbf{Proof}
For proving the necessary condition, observe that if $a_n = 0$ for $n \geq n_0$, then the series is a finite sum, hence $\theta$ is rational. If $a_n > 0$ infinitely often, let us to construct a nested sequence of closed intervals $I_n$ with intersection $\theta$. Let $I_1 = [a_0 + \frac{a_1}{b_1}, a_0 + \frac{a_1 + 1}{b_1}]$. Proceeding inductively, we have two possibilities, the first one, if $a_{n} = 0$, so define $I_{n} = I_{n-1}$. When $a_{n} \neq 0$, divide the interval $I_{n-1}$ into $b_{n} - a_{n} + 1$ ($\geq 2$) subintervals, the first one with length $\frac{a_{n}}{b_1\cdots b_{n}}$ and the other ones with equal length, namely, $\frac{1}{b_1\cdots b_{n}}$, and let the first one be $I_{n}$. By construction, $|I_n| \geq \frac{1}{b_1 \cdots b_{n}}$, for all $n \mem{N}$ and when $a_n \neq 0$, the length of $I_n$ is exactly $\frac{1}{b_1\cdots b_n}$. By hypothesis on $a_n$, there exist infinitely many $n \mem{N}$, such that $|I_n| = \frac{1}{b_1 \cdots b_n}$. Thus, we have
\begin{center}
$I_n = \left[a_0 + \frac{a_1}{b_1} + \ldots + \frac{a_n}{b_1 \cdots b_n}, a_0 + \frac{a_1}{b_1} + \ldots + \frac{a_n + 1}{b_1 \cdots b_n}\right] = \left[\frac{A_n}{b_1 \cdots b_n}, \frac{A_n + 1}{b_1 \cdots b_n}\right]$
\end{center}
where $A_n \mem{Z}$ for each $n \mem{N}$. Also $\theta \in I_n$ for all $n \geq 1$. In fact, by hypothesis it is easy see that $\theta > \frac{A_n}{b_1 \cdots b_n}$, for all $n \geq 1$. For the other inequality, note that  $\frac{a_{m}}{b_m} \leq 1 - \frac{1}{b_m}$, for all $m \mem{N}$, therefore 
\begin{equation}
b_1 \cdots b_n (\theta - (a_0 + \frac{a_1}{b_1} + \ldots + \frac{a_n}{b_1 \cdots b_n})) \leq 1
\end{equation}
Now if $a_n = b_n - 1$ for $n \geq n_0$, then $\theta$ is the right-hand endpoint of $I_{n_0 - 1}$, because each $I_n$ contains that endpoint and the lengths of the $I_n$ tend to zero. Hence again $\theta$ is rational. For showing the sufficient condition, note that if $a_m < b_m - 1$, then holds the strict inequality in $(2)$, for each $n < m$. Since $a_n > 0$ holds infinitely often,
\begin{center}
$\displaystyle\bigcap_{n = 1}^{\infty}I_n = \theta.$
\end{center}
Suppose that $\theta = \frac{p}{q} \mem{Q}$. Each prime number divides infinitely many $b_n$, so there exist $n_0$ sufficiently large such that $q | b_1 \cdots b_{n_0}$ and $a_{n_0} \neq 0$. Hence $b_1 \cdots b_{n_0} = kq$ for some $k \mem{N}$. Take $N \geq n_0$, such that, $a_{N+1} < b_{N+1} - 1$. Hence $\theta$ lies in interior of $I_{N}$. Also $I_N = I_{n_0 + k}$ for some $k \geq 0$. Suppose $I_N = I_{n_0}$. We can write $\theta = \frac{kp}{b_1 \cdots b_{n_0}}$, thus $\frac{A_{n_0}}{b_1 \cdots b_{n_0}} < \frac{kp}{b_1 \cdots b_{n_0}} < \frac{A_{n_0} + 1}{b_1 \cdots b_{n_0}}$. But that is a contradiction. If $I_N = I_{n_0+k}$, for $k \geq 1$, then we write $\theta = \frac{kpb_{n_0+1}\cdots b_{n_0+k}}{b_1 \cdots b_{n_0+k}}$. But that is again a contradiction. Therefore, it follows the irrationality of $\theta$.
\qed


\section{Irrationality measure}

The next step is to give an irrationality measure for some Cantor series. Now, we construct an uncountable family of functions, where one of them is exactly a well-known function for us.
\begin{defi}
Given $\sigma = (b_1, b_2, ...) \mem{N}^{\infty}$, satisfying
\begin{center}
$(*)$ For all $p$ prime number, the set $\{n \mem{N}\ |\ p|b_n\}$ is infinite.
\end{center}
We define the function $D(\cdot,\sigma):\mathbb{Z}^{*}\rightarrow \mathbb{N}$, by
\begin{center}
$D(q,\sigma):= \min\{n \mem{N}\ |\ q | b_1\cdots b_n\}$
\end{center}
Note that $D(\cdot,\sigma)$ is well defined, by condition $(*)$ and the well-ordering theorem. 
\end{defi}
In \cite{Sondow}, J. Sondow showed that for all integers $p$ and $q$ with $q > 1$, 
\begin{equation}\label{sondow}
\left|e - \displaystyle\frac{p}{q}\right| > \displaystyle\frac{1}{(S(q) + 1)!},
\end{equation}
where $S(q)$ is the smallest positive integer such that $S(q)!$ is a multiple of $q$ (the so-called Smarandache function, see \cite{Weisstein}). Note that if $\eta = (1, 2, 3, ...)$, then $D(q,\eta) = S(q)$. Since $e$ is a Cantor series and $D(\cdot, \sigma)$ is a generalization of Smarandache function, it is natural to think in a generalization or an improvement to the inequality in (3).

\textbf{Lemma} \textit{Given $n \in \mathbb{N}$, we have
\begin{equation}\label{desi1}
\left|\theta - \frac{m}{b_1 \cdots b_n}\right| \geq \min\left\{\left|\theta - \frac{A_n}{b_1 \cdots b_n}\right|, \left|\theta - \frac{A_n + 1}{b_1 \cdots b_n}\right|\right\}
\end{equation}
for all $m \mem{Z}.$}

\textbf{Proof}
Suppose that the result fail for some $m \mem{Z}$. So, $\frac{m}{b_1 \cdots b_n}$ lies in interior of $I_n$. Contradiction. Hence $(4)$ holds for all $m \mem{Z}$. 
\qed

\textbf{Proposition} \textit{Suppose that a Cantor series $\theta$, like in (1) and satisfying $(*)$, is an irrational number. For all integers $p\mem{Z}$ and $q \mem{Z}^*$, with $D(q, \sigma) > 1$, let $k$ be the smallest integer greater than $D(q,\sigma)$ such that the interval $I_k$ lies in the interior of $I_{D(q,\sigma)}$. Then 
\begin{equation}
\left |\theta - \displaystyle\frac{p}{q} \right | > \displaystyle\frac{\min\{a_k, b_k - a_k - 1\}}{b_1\cdots b_k}
\end{equation}
where $\sigma = (b_1, b_2, ...)$.}

\textbf{Proof}
Let $\sigma = (b_1, b_2, ...)$. Set $n = D(q, \sigma)$ and $m = \frac{pb_1\cdots b_n}{q}$. Therefore $m$ and $n$ are integers and 

\begin{eqnarray}
\left|\theta - \displaystyle\frac{p}{q} \right|  & = &
\left|\theta - \displaystyle\frac{m}{b_1 \cdots b_n} \right|\nonumber\\
 & \geq &
 \min\left\{\left|\theta - \frac{A_n}{b_1 \cdots b_n}\right|, \left|\theta - \frac{A_n + 1}{b_1 \cdots b_n}\right|\right\}\\
 & > &
\displaystyle\frac{\min\{a_k, b_k - a_k - 1\}}{b_1\cdots b_k}.\\
\nonumber
\end{eqnarray}
The inequalities $(6)$ and $(7)$ follow respectively by Lemma 1 and the hypothesis on $k$.

\qed 

The result below gives a slight improvement to (3).

\textbf{Corollary} \textit{If $p$ and $q$ are integers, with $q \neq 0$, then
\begin{equation}
\left|e - \displaystyle\frac{p}{q}\right| > \displaystyle\frac{1}{(D(q, \sigma) + 2)!},
\end{equation}
where $\sigma = (2, 3, 4, ...)$.}

\textbf{Proof}
Since that $\min_{p \mem{Z}}|e - p| > 0.28 > \frac{1}{6}$, then $(8)$ holds in the case $q = \pm 1$. In case $q \neq \pm 1$ the inequality also holds by Proposition and Example 1. Moreover, in this case we have $S(q) - 1 \in \{n \mem{N}\ | \ q | (n+1)!\}$ and $D(q, \sigma) + 1 \in \{n \mem{N}\ |\ q | n!\}$. Thus $S(q) = D(q,\sigma) + 1$. Hence
\begin{center}
$\left|e - \displaystyle\frac{p}{q}\right| > \displaystyle\frac{1}{(D(q, \sigma) + 2)!} = \displaystyle\frac{1}{(S(q) + 1)!}$
\end{center}
\qed

Actually, the improvement happens only because $(8)$ also holds for $q = \pm 1$.

\begin{exe}
The number $\xi := \frac{1}{(1!)^5} + \frac{1}{(2!)^5} + \frac{1}{(3!)^5} + \ldots = 1.031378...$ is irrational, moreover for $p,q \mem{Z}$ with $q \neq 0$, we have
\begin{center}
$\left|\xi - \displaystyle\frac{p}{q} \right| > \displaystyle\frac{1}{(D(q, \sigma) + 2)!^5}$
\end{center}
where $\sigma = (2^5, 3^5, ...)$.
\end{exe}

\textbf{Acknowledgments}

The author would like to thank Jonathan Sondow, Luiz Ant\^ onio Monte and Ana Paula Chaves for their helpful comments. The author is supported by CNPq.




\end{document}